\documentclass[reqno, 12pt, reqno]{amsart}

\usepackage{amsfonts, amsthm, amsmath, amssymb}
\usepackage{hyperref}
\hypersetup{colorlinks=false}

\usepackage[margin=1.5in]{geometry}

\usepackage{helvet}

\RequirePackage{mathrsfs} \let\mathcal\mathscr

\newtheorem{theorem}{Theorem}

\theoremstyle{definition}
\newtheorem*{ack}{Acknowledgement}

\renewcommand{\phi}{\varphi}
\renewcommand{\rho}{\varrho}

\renewcommand{\leq}{\leqslant}

\renewcommand{\bar}{\overline}

\title[A note on Burgess bound]{A note on Burgess bound}

\author{Ritabrata\ Munshi}
\address{School of Mathematics\\ 
Tata Institute of Fundamental Research\\
1 Homi Bhabha Road\\ Colaba\\Mumbai 400005\\ India}
\curraddr{Statistics and Mathematics Unit\\ Indian Statistical Institute\\ 203 B.T. Road\\ Kolkata 700108\\ India} 
\email{rmunshi@math.tifr.res.in}
\thanks{The author is supported by SwarnaJayanti Fellowship, 2011-12, DST, Government of India.}


\begin{document}

\begin{abstract}
Let $f$ be a $SL(2,\mathbb Z)$ Hecke cusp form, and let $\chi$ be a primitive Dirichlet character modulo $M$, which we assume to be prime. We prove the Burgess type bound for the twisted $L$-function:
\begin{align*}
L\left(\tfrac{1}{2},f\otimes\chi\right)\ll_{f,\varepsilon} M^{1/2-1/8+\varepsilon}.
\end{align*}
The method also yields the original bound of Burgess for Dirichlet $L$-functions:
\begin{align*}
L\left(\tfrac{1}{2},\chi\right)\ll_{\varepsilon} M^{1/4-1/16+\varepsilon}.
\end{align*}
\end{abstract}

\maketitle


\section{Introduction}
\label{introd}

In the series of papers \cite{Mu}, \cite{Mu4}, \cite{Mu0}, \cite{Mu5} and \cite{Mu16} a new approach to prove subconvexity has been proposed. This method has turned out to be quite effective in the case of degree three $L$-functions. Indeed if $\pi$ is a $SL(3,\mathbb{Z})$ Hecke-Maass cusp form then the $t$ aspect subconvex bound
\begin{align*}
L(\tfrac{1}{2}+it,\pi)\ll (2+|t|)^{3/4-1/16+\varepsilon},
\end{align*}
was established in \cite{Mu0},
and if $\chi$ is a primitive Dirichlet character modulo a prime $M$, then the twist aspect subconvex bound
\begin{align*}
L(\tfrac{1}{2},\pi\otimes \chi)\ll M^{3/4-1/308+\varepsilon},
\end{align*}
was established in \cite{Mu5} and \cite{Mu16}.
In the special case, where $\pi$ is a symmetric square lift of a $SL(2,\mathbb{Z})$ form, then the former bound was previously obtained by Li \cite{L}, building on the work of Conrey and Iwaniec \cite{CI}. Again for $\pi$ a symmetric square lift and $\chi$ a quadratic character, a subconvex bound (with a much stronger exponent) in the latter case was previously obtained by Blomer \cite{B}. The works of Li and Blomer fail to extend to the general case as the non-negativity of the central value of $L$-functions, which plays a vital role in their approach, ceases to hold in general. In contrast the method of \cite{Mu0}, \cite{Mu5}, \cite{Mu16} does not rely on this delicate feature. \\

Recently, Aggarwal \cite{Agar} showed that the computations in \cite{Mu0} when adapted in the $GL(2)$ scenario yields the Burgess exponent in the $t$ aspect. Taking a step further, refining certain estimates for exponential integral in \cite{Mu0}, Singh \cite{Sin} established the Weyl exponent in the $t$ aspect for $GL(2)$ $L$-functions 
\begin{align*}
L(\tfrac{1}{2}+it,g)\ll (2+|t|)^{1/3+\varepsilon}.
\end{align*}
Here $g$ is a Hecke modular or a Maass cuspform.
This is of course the celebrated result of Good \cite{Go}, which still remains the record though the corresponding exponent for the Riemann zeta function $\zeta(s)$ has gone through several improvements. It is now natural to ask what would one get if one had adopted the method of \cite{Mu5}, \cite{Mu16}, in the $GL(2)$ context. The aim of the present note is to give a brief analysis of this possibility. We show that the method yields the Burgess type bound.\\

\begin{theorem}
\label{mthm}
Suppose $g$ is a $SL(2,\mathbb{Z})$ Hecke cusp form, and $\chi$ is a primitive Dirichlet character modulo $M$ (which we assume to be prime). Then we have
\begin{align}
\label{main-bound}
L\left(\tfrac{1}{2},g\otimes\chi\right)\ll M^{1/2-1/8+\varepsilon}.
\end{align}
\end{theorem}
  
\bigskip

The exponent $1/2+\varepsilon$ corresponds to the convexity bound, which follows easily from the functional equation and basic complex analysis. The Riemann Hypothesis for $L(s,g\otimes\chi)$ implies the much stronger bound - the Lindelof Hypothesis - 
\begin{align*}
L\left(\tfrac{1}{2},g\otimes\chi\right)\ll M^{\varepsilon}.
\end{align*}
Our result \eqref{main-bound} is of course not new. Blomer and Harcos \cite{BH} established the Burgess exponent using Bykovskii's method. The weaker exponent $1/2-(1-2\theta)/8$, where $\theta$ is the exponent towards the Ramanujan conjecture, was previously established by Blomer, Harcos and Michel \cite{BHM}. The latter exponent becomes Burgess under the Ramanujan conjecture $\theta=0$. The first non-trivial bound of this sort, albeit with a weaker exponent, was obtained by Duke, Friedlander and Iwaniec \cite{DFI-1}. The Burgess bound, even in the original set up of Burgess \cite{Bu} for Dirichlet $L$-functions, still remains the record (though it has been improved in certain special cases - see e.g. \cite{BMi}, \cite{CI}, \cite{Mil}, \cite{MS}). Any improvement over this will be a major breakthrough in the field.\\

Our analysis should work even when $M$ is not a prime, but for simplicity we only provide a sketch for $M$ prime. From the technical point of view this paper puts the $GL(2)$ delta method of \cite{Mu6}, \cite{Mu16} in a better perspective. We closely follow \cite{Mu16}. In particular, we use the $GL(2)$ analogue of the `congruence-equation trick' (see \cite{BM}, \cite{Mu4}, \cite{Mu6}) and the trick of splitting  a dual variable of summation to optimize the saving after an application of Cauchy inequality. The $GL(2)$ set up simplifies certain technicalities, e.g. there is no need to sum over the modulus after dualization. However the extra sum over the level of forms (i.e. $p$ below) becomes redundant. This phenomenon is typical in $GL(2)$ set up, and that this sum provides the extra harmonics needed in the degree three case indicates that one should try to understand this in the amplification technique set up so as to extend that technique beyond $GL(2)$. We also observe that the cuspidality of $g$ is not required in the proof, and the proof holds even when the Fourier coefficient $\lambda_g(n)$ is replaced by the divisor function $d(n)$. In this way we can derive the classic result of Burgess.\\

\begin{theorem}
\label{mthm-2}
Suppose $\chi$ is a primitive Dirichlet character modulo $M$ (which we assume to be prime). Then we have
\begin{align}
\label{main-bound-2}
L\left(\tfrac{1}{2},\chi\right)\ll M^{1/4-1/16+\varepsilon}.
\end{align}
\end{theorem}

\bigskip

\ack
The author wishes to thank Roman Holowinsky for several insightful comments. He thanks the referee for helpful suggestions.\\


\section{Idea of the proof}

To get a subconvex bound for $L(1/2,\chi)^2$ one is led via the approximate functional equation to consider sums of the form
\begin{align*}
	S(N)=\sum_{n\sim N}d(n)\chi(n)
\end{align*}
where $d(n)$ is the divisor function and $N\ll M^{1+\varepsilon}$. We want to establish a bound of the form $|S(N)|/\sqrt{N}\ll M^{1/2-\theta}$ for some $\theta>0$. Roughly speaking, our first step is to rewrite the sum $S(N)$ as
\begin{align*}
	L^{-1}\sum_{\ell\in \mathcal{L}}\:\sum_{n\sim NL}\:\sum_{r\sim N}d(n)\chi(r)\delta(n,r\ell)
\end{align*}
where $\mathcal{L}$ is a set of $L$ primes of size $L^{1+\varepsilon}$, and $\delta(.,.)$ is the Kronecker delta symbol. Our job is to save $N^{3/2}LM^{\theta-1/2}$ for some $\theta>0$. To this end we use the harmonics from the space $S_k(pM,\psi)$ of cusp forms of weight $k$, level $pM$ and nebentypus $\psi$, to detect the equation $n=r\ell$ (see \eqref{circ-meth}). Here $p$ is a fixed prime and $\psi$ is a non-primitive odd character modulo $pM$ of conductor $p$. This yields two terms - the off-diagonal contribution involving Kloosterman sums and the dual contribution involving the Fourier coefficients of cusp forms. The off-diagonal contribution is negligibly small if we pick $p\gg NLM^\varepsilon/M$, as it involves the $J$-Bessel function 
\begin{align*}
	J_{k-1}\left(\frac{4\pi\sqrt{nr\ell}}{cpM}\right),
\end{align*}
with $c$ a positive integer. Note that we are taking $k$ to be large, like $1/\varepsilon$, and one has the bound $J_{k-1}(x)\ll x^{k-1}$.\\

The dual term is given by
\begin{align*}
	\sum_{\psi\bmod{p}}\:\sum_{\ell\in \mathcal{L}}\:\sum_{f\in H_k(pM,\psi)}\:\omega_f^{-1}\sum_{n\sim NL}\:d(n)\lambda_f(n)\sum_{r\sim N}\overline{\lambda_f(r\ell)}\chi(r),
\end{align*}
where our job is to save $N^{3/2}LM^{\theta-1/2}$ for some $\theta>0$.
Next we apply summation formulas on the sum over $n$ and $r$. These can be derived, for example, from the functional equations of $L(s,f)^2$ and $L(s,f\otimes\chi)$ respectively. With this the above sum reduces to
\begin{align*}
	\sum_{\psi\bmod{p}}\:\sum_{\ell\in \mathcal{L}}\:\sum_{f\in H_k(pM,\psi)}\:\omega_f^{-1}\:g_\psi\sum_{n\sim p^2M^2/NL}\:d(n)\overline{\lambda_f(np\ell)}\sum_{r\sim pM^2/N}\lambda_f(r)\bar{\chi(r)},
\end{align*}
and we make a saving of size $(NL/pM)(N/\sqrt{p}M)$. Here $g_\psi$ stands for the Gauss sum associated with $\psi$. It now remains to save $p^{3/2}M^{3/2+\theta}/\sqrt{N}$. We now apply the Petersson formula. The diagonal is easily seen to be small due to size of the variables. The off-diagonal is roughly of the form
\begin{align*}
	\sum_{\psi\bmod{p}}\:g_\psi\sum_{\ell\in \mathcal{L}}\:\sum_{n\sim p^2M^2/NL}\:d(n)\sum_{r\sim pM^2/N}\overline{\chi(r)}\sum_{c\ll pM/N}S_\psi(np\ell,r;cpM).
\end{align*}
In the off-diagonal the Petersson formula saves $\sqrt{pM}/\sqrt{\text{size of}\:c}=\sqrt{N}$, and then applying the Poisson summation formula on the sum over $r$ we save $(pM^2/N)/\sqrt{p^2M^2/N}=M/\sqrt{N}$. The sum over $\psi$ saves $\sqrt{p}$ more, and with this we arrive at the expression
\begin{align*}
	\sum_{\ell\in \mathcal{L}}\:\sum_{n\sim p^2M^2/NL}\:d(n)\sum_{r\sim p}\:\sum_{c\ll pM/N}\:e\left(\frac{\bar{c}n\ell}{rM}\right)\bar{\chi(c)}\:\mathfrak{C}
\end{align*}
where
\begin{align*}
	\mathfrak{C}=\sum_{a\bmod{M}}\chi(\bar{a}+r)e\left(\frac{a\bar{c}n\ell}{M}\right).
\end{align*}
Our job now is to save $pM^{1/2+\theta}/\sqrt{N}$ in the above sum (beyond square root cancellation in the character sum $\mathfrak{C}$). Applying Cauchy inequality we seek to save $p^2M^{1+2\theta}/N$ in
\begin{align*}
	\sum_{n\sim p^2M^2/NL}\:\Bigl|\sum_{\ell\in \mathcal{L}}\:\sum_{r\sim p}\:\sum_{c\ll pM/N}\:e\left(\frac{\bar{c}n\ell}{rM}\right)\bar{\chi(c)}\:\mathfrak{C}\Bigr|^2.
\end{align*}
We will now open the absolute value and apply the Poisson summation formula on the sum over $n$. The diagonal contribution is seen to be satisfactory as long as we have enough terms inside the absolute value, namely $p^2LM/N>p^2M^{1+2\theta}/N$, i.e. $L>M^{2\theta}$. On the other hand the off-diagonal is satisfactory as long as $p^2M^2/NLM^{1/2}>p^2M^{1+2\theta}/N$, i.e. $L<M^{1/2-2\theta}$. At this stage we encounter a complete character sum of the form 
\begin{align*}
	\sum_{x\in \mathbb{F}_M-\{\text{few points}\}}\:\chi\left(\frac{P(x)}{Q(x)}\right)
\end{align*}
where $P$ and $Q$ are quadratic polynomials. Exactly the same sum appeared in Burgess' method, and like him we appeal to Weil's results (Riemann hypothesis for curves over finite fields) to conclude square-root cancellation in the sum. Now we observe that the optimal choice for $L$ is given by $L=M^{1/4}$ and $\theta$ is taken to be $1/8$. This establishes the bound
\begin{align*}
	\frac{1}{\sqrt{N}}\left|\sum_{n\sim N} d(n)\chi(n)\right|\ll M^{1/2-1/8+\varepsilon}
\end{align*}
for all $N$. From this we are able to conclude the Burgess bound
\begin{align*}
	L(1/2,\chi)\ll M^{3/16+\varepsilon}.
\end{align*}
A careful reader will observe that the above sketch works even when the divisor function is replaced by Fourier coefficients for of $GL(2)$ forms. Below we gives details for this more general case.\\


\section{The set up}

Let $g$ be a fixed cusp form for $SL(2,\mathbb{Z})$ with normalized Fourier coefficients $\lambda_g(n)$. Let $\chi$ be a primitive Dirichlet character modulo a prime $M$. The approximate functional equation for the twisted $L$-function yields 
\begin{align*}
L\left(\tfrac{1}{2},g\otimes \chi\right)\ll M^\varepsilon \sup_N \frac{|S(N)|}{N^{1/2}} + M^{(1-\theta)/2}
\end{align*}
where the supremum is taken over $M^{1-\theta}<N<M^{1+\varepsilon}$, and $S(N)$ are sums of the form
\begin{align*}
S(N)=\sum_{n=1}^\infty \lambda_g(n)\chi(n)W\left(\frac{n}{N}\right)
\end{align*}
with $W$ a smooth bump function supported on $[1,2]$. We will use the $GL(2)$ delta method to separate the oscillations of the Fourier coefficients and the character in the above sum.\\

We start by recalling the $GL(2)$ delta method. Let $p$ be a prime number and let $k\equiv 3\bmod{4}$ be a positive integer (which will be of the size $1/\varepsilon$). Let $\psi$ be a character of $\mathbb{F}_p^\times$ satisfying $\psi(-1)=-1=(-1)^k$. We consider $\psi$ as a character modulo $pM$. 
 The inclusion of $M$, which is the conductor of $\chi$, in the level is an analogue of the `congruence-equation trick'.  
Let $H^\star(pM,\psi)$ be the set of Hecke-normalized newforms and we extend it to $H_k(pM,\psi)$ - an orthogonal Hecke basis of the space of cusp forms $S_k(pM,\psi)$. 
Using the Petersson formula we derive
\begin{align}
\label{circ-meth}
&\delta(n,r)=\frac{1}{p-1}\;\sum_{\psi\bmod{p}}\left(1-\psi(-1)\right)\sum_{f\in H_k(pM,\psi)}\omega_f^{-1}\lambda_f(n)\overline{\lambda_f(r)}\\
\nonumber &-\frac{2\pi i}{p-1}\;\sum_{c=1}^\infty \frac{1}{cpM}\sum_{\psi\bmod{p}}\left(1-\psi(-1)\right)S_\psi(r,n;cpM)J_{k-1}\left(\frac{4\pi\sqrt{nr}}{cpM}\right),
\end{align}
i.e. the right hand side is $1$ if $n=r$, and is equal to $0$ otherwise. Here $\lambda_f(n)$ are the normalized Fourier coefficients of the form $f$, $\omega_f^{-1}$ is the spectral weight, 
$$
S_\psi(a,b;c)=\sum_{\alpha\bmod{c}}\psi(\alpha)e\left(\frac{\alpha a+\bar{\alpha}b}{c}\right)
$$ 
is the generalized Kloosterman sum and $J_{k-1}(x)$ is the Bessel function. (Note that we do not sum over $p$ as this fails to give any extra leverage.)\\

Let $\mathcal{L}$ be a set of primes in the range $[L,2L]$, such that
\begin{align*}
L^\star=\sum_{\ell\in\mathcal{L}}|\lambda_g(\ell)|^2\gg L^{1-\varepsilon},
\end{align*}
 and $L\ll M^{1-\varepsilon}$. Consider the sum
\begin{align}
\label{F}
\mathcal{F}=&\frac{1}{L^\star (p-1)}\;\sum_{\ell\in\mathcal{L}}\bar{\lambda_g(\ell)}\;\sum_{\psi\bmod{p}}\left(1-\psi(-1)\right)\sum_{f\in H_k(pM,\psi)}\omega_f^{-1}\\
\nonumber\times &\mathop{\sum\sum}_{m,n=1}^\infty \lambda_g(n)\lambda_f(n)\psi(m)W\left(\frac{nm^2}{N\ell}\right)\sum_{r=1}^\infty \overline{\lambda_f(r\ell)}\chi(r)V\left(\frac{r}{N}\right).
\end{align}
Here $W$ is a smooth bump function with support $[1,2]$ and $V$ is a smooth function supported in $[2/3,3]$, with $V(x)=1$ for $x\in [1,2]$, and satisfying $y^jV^{(j)}(y)\ll_j 1$. The diagonal term in the Petersson formula corresponds to $n=r\ell$, in which case we are left with the weight function $W(rm^2/N)V(r/N)$. This is clearly vanishing if $m\neq 1$. Hence the diagonal term is given by
\begin{align*}
&\frac{1}{L^\star}\;\sum_{\ell\in\mathcal{L}}\bar{\lambda_g(\ell)}\;\mathop{\sum}_{r=1}^\infty \lambda_g(r\ell)\chi(r)W\left(\frac{r}{N}\right)=S(N)+O\left(\frac{NM^\varepsilon}{L}\right).
\end{align*}
From \eqref{circ-meth} it follows that
\begin{align}
\label{s-star-n}
S(N)=\mathcal{F}-2\pi i\:\mathcal{O}+O\left(\frac{NM^\varepsilon}{L}\right)
\end{align}
where the off-diagonal is given by
\begin{align*}
\mathcal{O}=&\frac{1}{L^\star (p-1)}\;\sum_{\psi\bmod{p}}\left(1-\psi(-1)\right) \sum_{\ell\in\mathcal{L}}\bar{\lambda_g(\ell)}\:\mathop{\sum\sum}_{m,n=1}^\infty \lambda_g(n)\psi(m)W\left(\frac{nm^2}{N\ell}\right)\\
\nonumber &\times \sum_{r=1}^\infty \chi(r)V\left(\frac{r}{N}\right)\sum_{c=1}^\infty \frac{S_\psi(r\ell,n;cpM)}{cpM}J_{k-1}\left(\frac{4\pi\sqrt{n\ell r}}{cpM}\right).
\end{align*}
If we pick $p>LNM^\varepsilon/M$, then the off-diagonal becomes negligibly small due to the bound $J_{k-1}(x)\ll x^{k-1}$. So with this choice of $(p,L)$ our job is reduced to estimating the dual sum $\mathcal{F}$. 

\bigskip

  
\section{Treating the old forms}

To analyse the sum $\mathcal{F}$ we use  the functional equation for the $GL(2)\times GL(2)$ Rankin-Selberg convolution. There are two issues. First we need to split $\lambda_f(r\ell)$ using the Hecke relation and secondly we need to take care of the oldforms. To this end let $\mathcal{F}^\sharp$ be same as the expression in \eqref{F} with $\lambda_f(r\ell)$ replaced by $\lambda_f(r)\lambda_f(\ell)$. On using the Hecke relation in \eqref{F} we get two terms, one of them being $\mathcal{F}^\sharp$. To tackle the other term, consider
\begin{align}
\label{hr}
\sum_{f\in H_k(pM,\psi)}\omega_f^{-1}\;\sum_{\ell\in\mathcal{L}}\bar{\lambda_g(\ell)}&\mathop{\sum\sum}_{m,n=1}^\infty \lambda_g(n)\lambda_f(n)\psi(m) W\left(\frac{nm^2}{N\ell}\right)\\
\nonumber &\times \sum_{r=1}^\infty \overline{\lambda_f(r)}\chi(r\ell)V\left(\frac{r\ell}{N}\right).
\end{align}
To this we apply the Petersson formula. The diagonal is given by
\begin{align*}
\sum_{\ell\in\mathcal{L}}\bar{\lambda_g(\ell)}&\mathop{\sum\sum}_{m,n=1}^\infty \lambda_g(n)\chi(n\ell)\psi(m) W\left(\frac{nm^2}{N\ell}\right)V\left(\frac{n\ell}{N}\right).
\end{align*}
In the above sum $m$ has to be of size $L$. The sum over $\psi$ gives a saving of $\min\{L,p\}$. It turns out that the contribution of the above diagonal to $\mathcal{F}$ is bounded by $O(NM^\varepsilon/\min\{p,L\})$. The off-diagonal is given by
\begin{align*}
\sum_{c=1}^\infty &\sum_{\ell\in\mathcal{L}}\bar{\lambda_g(\ell)}\mathop{\sum\sum}_{m,n=1}^\infty \lambda_g(n)\psi(m) W\left(\frac{nm^2}{N\ell}\right)\\
&\times \sum_{r=1}^\infty \chi(r\ell)\frac{S_\psi(r,n;cpM)}{cpM}J_{k-1}\left(\frac{4\pi\sqrt{rn}}{cpM}\right)V\left(\frac{r\ell}{N}\right).
\end{align*}
This is negligibly small if $p>M^\varepsilon$ as $J_{k-1}(x)\ll x^{k-1}$. Hence we get
\begin{align*}
\mathcal{F}=\mathcal{F}^\sharp + O\left(\frac{NM^\varepsilon}{\min\{p,L\}}\right).
\end{align*}\\

Next we take into account the contribution of the old forms in $\mathcal{F}^\sharp$. For $h\in S_k(p,\psi)$, we set $h|_M(z)=M^{k/2}h(Mz)$ which lies in $S_k(pM,\psi)$. Define 
$
h^\star=h|_M-\left<h|_M,h\right>\left<h,h\right>^{-1}h.
$
Then $\{h,h^\star: h\in H_k(p,\psi)\}$ gives an orthogonal Hecke basis of the space of oldforms. For $f=h$, the functional equations yield the bound
\begin{align*}
\mathop{\sum\sum}_{m,n=1}^\infty \lambda_g(n)\lambda_f(n)\psi(m)W\left(\frac{nm^2}{N\ell}\right)\sum_{r=1}^\infty \overline{\lambda_f(r)}\chi(r)V\left(\frac{r}{N}\right)\ll p\;p^{1/2}M^{1+\varepsilon},
\end{align*}
as $g\otimes f$ has conductor $p^2$ and $f\otimes \chi$ has conductor $pM^2$. Now to treat $f=g^\star$, note that $\lambda_{g|_M}(r)=0$ unless $M|r$. But in this case we have  $\chi(r)=0$. Also $\left<g|_M,g|_M\right>=\left<g,g\right>$, hence by Bessel inequality we have $|\left<g|_M,g\right>\left<g,g\right>^{-1}|\leq 1$. Consequently for $f=g^\star$, the above sum reduces exactly to the case $f=g$ with a multiplier of size $1$. Hence the above bound also holds for $f=g^\star$. Consequently the contribution of the oldforms is bounded by $p^{3/2}M^\varepsilon$. Hence we get
\begin{align}
\label{f-2-fstar}
\mathcal{F}=\mathcal{F}^\star +O\left(p^{3/2}M^{\varepsilon}+\frac{NM^\varepsilon}{\min\{p,L\}}\right),
\end{align}
where
\begin{align}
\label{Fstar}
\mathcal{F}^\star=&\frac{1}{L^\star (p-1)}\:\sum_{\psi\bmod{p}}\left(1-\psi(-1)\right)\sum_{f\in H_k^\star(pM,\psi)}\omega_f^{-1}\;\sum_{\ell\in\mathcal{L}}\bar{\lambda_f(\ell)}\bar{\lambda_g(\ell)}\\
\nonumber\times &\mathop{\sum\sum}_{m,n=1}^\infty \lambda_g(n)\lambda_f(n)\psi(m)W\left(\frac{nm^2}{N\ell}\right)\sum_{r=1}^\infty \overline{\lambda_f(r)}\chi(r)V\left(\frac{r}{N}\right).
\end{align}
It follows that
\begin{align*}
S(N)=\mathcal{F}^\star +O\left(p^{3/2}M^{\varepsilon}+\frac{NM^\varepsilon}{\min\{p,L\}}\right),
\end{align*}
if $p>\max\{LNM^\varepsilon/M,M^\varepsilon\}$.

\bigskip


\section{Applying functional equations}
\label{fe}

Next we will apply the functional equation for $GL(2)\times GL(2)$ Rankin-Selberg convolution, and $GL(2)$ Voronoi summation formula to the sums over $(m,n)$ and $r$ respectively. Applying Mellin inversion formula we get
\begin{align}
\label{sum-1}
\mathop{\sum\sum}_{m,n=1}^\infty \lambda_g(n)\lambda_f(n)\psi(m)W\left(\frac{nm^2}{N\ell}\right)=\frac{1}{2\pi i}\int_{(2)}(N\ell)^s\tilde{W}(s)L(s,f\otimes g)\mathrm{d}s,
\end{align}
where $\tilde{W}$ stands for the Mellin transform of $W$. Applying the functional equation (see \cite{KMV}) we get
\begin{align*}
\frac{1}{2\pi i}\int_{(2)}(N\ell)^s\tilde{W}(s)\;\eta^2 \left(\frac{pM}{4\pi^2}\right)^{1-2s}\frac{\gamma(1-s)}{\gamma(s)} L(1-s,f\otimes \bar{g})\mathrm{d}s,
\end{align*}
where $\gamma(s)=\Gamma(s+(k-\kappa)/2)\Gamma(s+(k+\kappa)/2-1)$ and $\kappa$ is the weight of the form $g$. The sign of the functional equation 
\begin{align*}
\eta^2=\frac{g_\psi^2\psi^2(M)}{\lambda_g(M^2p^2)Mp}=\frac{g_\psi^2\psi(M)}{\lambda_g(p^2)p}
\end{align*}
where $g_\psi$ is the Gauss sum associated with the character $\psi$. Next we expand the $L$-function into a Dirichlet series and take dyadic subdivision. Since we are taking $k$ to be large (like $1/\varepsilon$), the poles of the gamma factor are away from $0$. Hence by shifting contours to the right or left we can show that the contribution of the terms from the blocks with $nm^2\notin [p^2M^{2-\varepsilon}/NL,p^2M^{2+\varepsilon}/NL]$ is negligibly small.  Hence the sum in \eqref{sum-1} essentially gets transformed into
\begin{align}
\label{sum-11}
\varepsilon_\psi^2\psi(M)\:\frac{N\ell}{pM}\:\mathop{\sum\sum}_{m,n=1}^\infty \lambda_g(n)\bar{\lambda_f(np^2)}\bar{\psi}(m)W\left(\frac{nm^2}{\tilde N}\right)
\end{align}
where $\varepsilon_\psi$ is the sign of the Gauss sum $g_\psi$ and $\tilde N\asymp p^2M^2/NL$. (In this paper the notation $A\asymp B$ means that $B/M^\varepsilon\ll A\ll BM^\varepsilon$, with implied constants depending on $\varepsilon$.)\\

To the other sum
\begin{align*}
\sum_{r=1}^\infty \overline{\lambda_f(r)}\chi(r)V\left(\frac{r}{N}\right)=g_{\bar{\chi}}^{-1}\sum_{a\bmod{M}}\bar{\chi}(a)\sum_{r=1}^\infty \overline{\lambda_f(r)}e\left(\frac{ar}{M}\right)V\left(\frac{r}{N}\right)
\end{align*}
we apply the Voronoi summation formula. This transforms the above sum into
\begin{align*}
2\pi i^k\chi(-p)\bar{\psi}(-M)\frac{g_\chi g_{\bar{\psi}}}{g_{\bar{\chi}}}\frac{N}{Mp}\sum_{r=1}^\infty \lambda_f(rp)\bar{\chi}(r)\int_0^\infty V(x)J_{k-1}\left(\frac{4\pi\sqrt{N rx}}{M\sqrt{p}}\right)\mathrm{d}x.
\end{align*}
As $k$ is large the Bessel function is negligibly small if $r\ll M^2P/NM^\varepsilon$. On the other hand making the change of variables $y^2=x$, pulling out the oscillation of the Bessel function and integrating by parts we get that the integral is negligibly small if $r\gg M^{2+\varepsilon}P/N$.
This reduces the analysis of the sum in \eqref{Fstar} to that of the sums of the type
\begin{align}
\label{to-ana}
\mathcal{D}^\star=&\frac{N^2}{M^2p^3}\: \sum_{\psi\bmod{p}}\left(1-\psi(-1)\right)g_{\psi}\sum_{f\in H_k^\star(pM,\psi)}\omega_f^{-1}\overline{\lambda_f(p)}\sum_{\ell\in\mathcal{L}}
\bar{\lambda_f(\ell)}\bar{\lambda_g(\ell)}\\
\nonumber\times &\mathop{\sum\sum}_{m,n=1}^\infty \lambda_g(n)\overline{\lambda_f(n)}\bar{\psi(m)} W\left(\frac{nm^2}{\tilde N}\right)\sum_{r=1}^\infty \lambda_f(r)\bar{\chi}(r)V\left(\frac{r}{\tilde{R}}\right)
\end{align}
where
\begin{align}
\label{range-all}
\frac{p^2M^2}{NLM^{\varepsilon}}\ll \tilde{N}\ll \frac{p^2M^{2+\varepsilon}}{NL}\;\;\;\text{and}\;\;\; \frac{pM^{2}}{NM^\varepsilon}\ll \tilde{R}\ll \frac{pM^{2+\varepsilon}}{N}.
\end{align}\\

Next we shall apply the Petersson formula, to this end we first need to extend the sum over $f$ to a full orthogonal basis. So we estimate the contribution of the oldforms.
If $f$ is an oldform coming from level $p$, the sub sum over $(m,n)$ in \eqref{to-ana}  is negligibly small, as the length of the sum $\tilde{N}\gg P^{2}M^\varepsilon$ is larger than the size of the conductor $P^2$. (Note that $L<M^{1-\varepsilon}$.) So the sum over $f$ can be extended to a complete Hecke basis at a cost of a negligible error term. Next we use the Hecke relation. We analyse the generic term. The other term can be analysed in the same fashion and at the end we get a stronger bound for it. Consider \eqref{to-ana} with $\lambda_f(p)\lambda_f(\ell)\lambda_f(n)$  replaced by $\lambda_f(p\ell n)$  and the sum over $f$ is extended to a full orthogonal Hecke basis - 
\begin{align}
\label{to-ana-dual-sum}
\mathcal{D}=&\frac{N^2}{M^2p^3}\: \sum_{\psi\bmod{p}}\left(1-\psi(-1)\right)g_{\psi}\sum_{f\in H_k(pM,\psi)}\omega_f^{-1}\sum_{\ell\in\mathcal{L}}\bar{\lambda_g(\ell)}\\
\nonumber\times &\mathop{\sum\sum}_{m,n=1}^\infty \lambda_g(n)\overline{\lambda_f(np\ell)}\bar{\psi(m)} W\left(\frac{nm^2}{\tilde N}\right)\sum_{r=1}^\infty \lambda_f(r)\bar{\chi}(r)V\left(\frac{r}{\tilde{R}}\right).
\end{align}
We apply the Petersson formula. The diagonal is given by 
\begin{align*}
\Delta=&\frac{N^2}{M^2p^3}\: \sum_{\psi\bmod{p}}\left(1-\psi(-1)\right)g_{\psi}\:\sum_{\ell\in\mathcal{L}}\bar{\lambda_g(\ell)}\\
\nonumber\times &\mathop{\sum\sum}_{m,n=1}^\infty \lambda_g(n)\bar{\psi(m)}\bar{\chi(np\ell)} W\left(\frac{n m^2}{\tilde N}\right)V\left(\frac{np\ell}{\tilde{R}}\right).
\end{align*}
The sum over $\psi$ can now give a saving of $p^{1/2}$, and so it follows that 
\begin{align*}
\Delta\: \ll \;\frac{NM^\varepsilon}{p}.
\end{align*}
The off-diagonal is given by
\begin{align}
\label{123d}
\mathcal{O}^\star=&\frac{N^2}{M^2p^3}\: \sum_{\psi\bmod{p}}\left(1-\psi(-1)\right)g_{\psi}\:\sum_{\ell\in\mathcal{L}}\bar{\lambda_g(\ell)}\;\sum_{r=1}^\infty \bar{\chi}(r)V\left(\frac{r}{\tilde{R}}\right)\\
\nonumber\times &\mathop{\sum\sum}_{m,n=1}^\infty \lambda_g(n)\bar{\psi(m)} W\left(\frac{nm^2}{\tilde N}\right)\sum_{c=1}^\infty \frac{S_\psi(np\ell,r;cpM)}{cpM}J_{k-1}\left(\frac{4\pi\sqrt{nr\ell}}{cM\sqrt{p}}\right).
\end{align}
Using the bound for the $J$-Bessel function, we can now cut the $c$ sum at $c\ll C=pM^{1+\varepsilon}/Nm$ at the cost of a negligible error term.\\

\section{Applying Poisson summation}
\label{poisson}

Next we apply the Poisson summation formula on the sum over $r$. This yields
\begin{align*}
\mathcal{O}^\star=&\frac{N^2}{M^2p^3}\: \sum_{\ell\in\mathcal{L}}\bar{\lambda_g(\ell)}\;\mathop{\sum\sum}_{m,n=1}^\infty \lambda_g(n) W\left(\frac{nm^2}{\tilde N}\right)\:\sum_{c=1}^\infty \frac{\tilde{R}}{(cpM)^2}\sum_{r\in \mathbb{Z}} 
\;\mathfrak{C}\:\mathfrak{I}
\end{align*}
where the character sum is given by
\begin{align*}
\mathfrak{C}=\sum_{\psi\bmod{p}}\left(1-\psi(-1)\right)g_{\psi}\:\bar{\psi}(m)\;\sum_{a\bmod{cpM}}\bar{\chi}(a)S_\psi(np\ell,a;cpM)e\left(\frac{ar}{cpM}\right),
\end{align*}
and the integral is given by
\begin{align*}
\mathfrak{I}=\int V\left(y\right)J_{k-1}\left(\frac{4\pi\sqrt{n\tilde{R}\ell y}}{cM\sqrt{p}}\right)e\left(-\frac{\tilde{R}ry}{cpM}\right)\mathrm{d}y.
\end{align*}
Since $(M,cp)=1$ the character sum splits as a product of two character sums - one modulo $M$ and the other modulo $cp$. In the part modulo $cp$, 
\begin{align*}
\sum_{\psi\bmod{p}}g_{\psi}\:\bar{\psi}(\pm m)\;\sum_{a\bmod{cp}}S_\psi(\bar{M} np\ell,\bar{M} a;cp)e\left(\frac{ar\bar{M}}{cp}\right),
\end{align*}
after opening the Kloosterman sum, we get
\begin{align*}
\sum_{\psi\bmod{p}}g_{\psi}\:\bar{\psi}(\pm m)\;\sideset{}{^\star}\sum_{b\bmod{cp}}\psi(b)e\left(\frac{b\bar{M}np\ell}{cp}\right)\sum_{a\bmod{cp}}e\left(\frac{\bar{bM}a}{cp}+\frac{ar\bar{M}}{cp}\right).
\end{align*}
The innermost sum vanishes unless $(r,cp)=1$ (which rules out the existence of zero frequency $r=0$), in which case the above sum reduces to
\begin{align*}
cp\:\sum_{\psi\bmod{p}}g_{\psi}\:\bar{\psi}(\mp mr)\;e\left(-\frac{\bar{rM}np\ell}{cp}\right).
\end{align*}
Evaluating the $\psi$ sum we arrive at $\mathfrak{C}=\mathfrak{C}_{-}-\mathfrak{C}_+$ where
\begin{align*}
\mathfrak{C}_\pm = cp(p-1)\:e\left(\pm\frac{mr}{p}-\frac{\bar{rM}n\ell}{c}\right)\:\sum_{a\bmod{M}}\bar{\chi}(a)S(\bar{c} n\ell,\bar{cp} a;M)e\left(\frac{ar\bar{cp}}{M}\right).
\end{align*}
To analyse the integral $\mathfrak{I}$ we extract the oscillation from the Bessel function, this transforms the integral to 
\begin{align*}
\int W\left(y\right)\:e\left(\pm\frac{2\sqrt{n\tilde{R}\ell y}}{cM\sqrt{p}}-\frac{\tilde{R}ry}{cpM}\right)\mathrm{d}y,
\end{align*}
where $W$ is compactly supported and satisfies $y^jW^{(j)}(y)\ll 1$. The stationary phase is given by $y_0=n\ell p/\tilde{R}r^2$. Expanding at the stationary phase we are able to replace the integral $\mathfrak{I}$ by
\begin{align*}
\frac{cM\sqrt{p}}{\sqrt{n\tilde{R}\ell}}\: e\left(\frac{n\ell}{cMr}\right)\:W(\dots)
\end{align*}
where the weight function $W$ is non-oscillatory. Moreover $W$ is negligibly small unless $|r|\asymp p/m$, and we can separate the variables in $W$ at the cost of a negligible error term. (Note that it follows that the contribution of the terms with $m\gg pM^\varepsilon$ is negligibly small.)  Then applying the reciprocity relation
\begin{align*}
e\left(-\frac{\bar{rM}n\ell}{c}\right)e\left(\frac{n\ell}{cMr}\right)=e\left(\frac{\bar{c}n\ell}{rM}\right)
\end{align*}
 we see that we have transformed the off-diagonal $\mathcal{O}^\star$ to 
\begin{align*}
\mathcal{O}^\dagger=&\frac{N^2}{M^3p^3}\: \sum_{\ell\in\mathcal{L}}\bar{\lambda_g(\ell)}\;\mathop{\sum\sum}_{nm^2\sim \tilde{N}} \lambda_g(n)\:m\:\sum_{c\ll C} \:\sum_{\substack{|r|\asymp p/m\\(r,cp)=1}} 
e\left(\pm\frac{mr}{p}+\frac{\bar{c}n\ell}{rM}\right)\;\mathcal{C}
\end{align*}
where
\begin{align*}
\mathcal{C}&=\sum_{a\bmod{M}}\bar{\chi}(a)S(\bar{c} n\ell,\bar{c} a;M)e\left(\frac{ar\bar{c}}{M}\right)\\
&=g_{\bar{\chi}}\:\bar{\chi(c)}\:\sideset{}{^\star}\sum_{\alpha\bmod{M}}\chi(\bar{\alpha}+r)e\left(\frac{\alpha\bar{c}n\ell}{M}\right),
\end{align*}
and the weight function is non-oscillatory.
From the analysis in this section we now conclude that
\begin{align*}
S(N)\ll |\mathcal{O}^\dagger| +p^{3/2}M^{\varepsilon}+\frac{NM^\varepsilon}{\min\{p,L\}},
\end{align*}
if $p>\max\{LNM^\varepsilon/M,M^\varepsilon\}$, and $L<M^{1-\varepsilon}$.

\bigskip

\section{Applying Cauchy inequality followed by Poisson summation}

Applying the Cauchy inequality we get
\begin{align*}
\mathcal{O}^\dagger\ll \frac{N^2\tilde{N}^{1/2}}{M^3p^3}\;\sum_{m\ll\tilde{N}}\:\Omega^{1/2}
\end{align*}
where
\begin{align*}
\Omega=&\mathop{\sum}_{n\sim \tilde{N}/m^2}\:\left|\sum_{\ell\in\mathcal{L}} \bar{\lambda_g(\ell)}\:\sum_{c\ll C} \:\sum_{\substack{|r|\asymp p/m\\(r,cp)=1}} 
e\left(\pm\frac{mr}{p}+\frac{\bar{c}n\ell}{rM}\right)\;\mathcal{C}\right|^2.
\end{align*}
We extend the outer sum to all integers and insert a suitable bump function $W(nm^2/\tilde{N})$. Then opening the absolute value square and applying Poisson summation on the sum over $n$ we arrive at
\begin{align*}
\Omega\ll & \frac{\tilde{N}}{m^2}\:\mathop{\mathop{\sum\sum}_{\ell_1,\ell_2\in \mathcal{L}}\:\mathop{\sum\sum}_{c_1,c_2\ll C} \:\mathop{\sum\sum}_{|r_1|,|r_2|\asymp p/m}\; \sum_{|n|\ll p^2M^{1+\varepsilon}/\tilde{N}}}_{ F(\dots)\equiv 0\bmod{r_1r_2}}
\;|\mathcal{C}^\dagger|.
\end{align*}
where $F=F(\dots)=\ell_1c_2r_2-\ell_2c_1r_1+c_1c_2n$ and 
\begin{align*}
\mathcal{C}^\dagger =\sum_{b\bmod{M}}&\mathop{\sum\sum}_{\substack{a_1,a_2\bmod{M}\\(a_1a_2,M)=1}}\chi(\bar{a_1}+r_1)\bar{\chi}(\bar{a_2}+r_2)
e\left(\frac{a_1b\bar{c_1}\ell_1}{M}-\frac{a_2b\bar{c_2}\ell_2}{M}\right)\\
&\times e\left(\frac{\bar{c_1}\ell_1b\bar{r_1}}{M}-\frac{\bar{c_2}\ell_2b\bar{r_2}}{M}+\frac{b\bar{r_1r_2}n}{M}\right).
\end{align*}
Evaluating the sum over $b$ the character sum reduces to
\begin{align*}
M\sum_{x\in \mathbb{F}_M-\{\dots\}}\chi(\bar{c_1}\ell_1\overline{(x\bar{c_2}\ell_2-\bar{c_1c_2r_1r_2}F)}+r_1)\bar{\chi}(\bar{x}+r_2).
\end{align*}
The sum over $x$ misses a bounded number of points, and in fact can be rewritten as
\begin{align*}
\sum_{\substack{x\in \mathbb{F}_M-\{\dots\}}}\chi\left(\frac{Q_1(x)}{Q_2(x)}\right)
\end{align*}
where $Q_i$ are quadratic polynomials. The Riemann hypothesis for curves (Weil's theorem) gives square-root cancellation in such type of sums in general.\\

We will first estimate the contribution of the terms for which $F\equiv 0\bmod{M}$. In this case the character sum reduces further to 
\begin{align*}
M\sum_{x\in \mathbb{F}_M-\{\text{few points}\}}\chi(\bar{c_1}\ell_1c_2\bar{\ell_2}\bar{x}+r_1)\bar{\chi}(\bar{x}+r_2).
\end{align*}
This is bounded by $O(M)$ unless $r_1=r_2$ and $c_1\ell_2=c_2\ell_1$ in which case it is bounded by $O(M^2)$. (We are assuming that $pL<NM^{-\varepsilon}$.) As $(r_1r_2,M)=1$ we have $F\equiv 0 \bmod{r_1r_2M}$, but since here the modulus is larger than the size of the equation (which is $p^2ML/Nm^2$), we get $F=0$. Now we need to count the number of solutions of $F=0$. If we have the restriction $r_1=r_2$ and $c_1\ell_2=c_2\ell_1$, then $F=0$ implies that $n=0$ (the zero frequency). So the total contribution of this case to $\Omega$ is bounded by $O(\tilde{N}M^{2+\varepsilon}LCp/m^2)$. Otherwise $F=0$ reduces to the congruence 
$$
\ell_1c_2r_2-\ell_2c_1r_1\equiv 0\bmod{c_1c_2}.
$$
If we define $d=(c_1,c_2)$ and write $c_i=dc_i'$, then $c_1'|\ell_1r_2$, $c_2'|\ell_2r_1$
Hence the contribution of this case to $\Omega$ is dominated by
\begin{align*}
\frac{\tilde{N}M}{m^2}\:\sum_{d\ll C}\mathop{\sum\sum}_{\ell_1,\ell_2\in \mathcal{L}} \:\mathop{\sum\sum}_{|r_1|,|r_2|\asymp p/m}
\;1 \ll \frac{\tilde{N}ML^2p^2C}{m^2}.
\end{align*}
The last bound is smaller than the bound we obtained above for the zero frequency as $pL<M$. \\

Now consider the case where $M\nmid F$. Here Weil bound gives $\mathcal{C}^\dagger\ll M^{3/2}$ (see Corollary~11.24 of \cite{IK}). In this case the contribution of $n=0$ to $\Omega$ is bounded by 
\begin{align*}
\frac{\tilde{N}M^{3/2}}{m^2}\:\sum_{d\ll p/m}\mathop{\sum\sum}_{\ell_1,\ell_2\in \mathcal{L}} \:\mathop{\sum\sum}_{|c_1|,|c_2|\ll C}
\;1 \ll \frac{\tilde{N}M^{3/2}L^2pC^2}{m^2}.
\end{align*}
Now suppose $n\neq 0$. We will use the congruence condition to count the number of $\ell_i$. Indeed the congruence condition splits into two conditions 
\begin{align*}
\ell_2c_1r_1\equiv c_1c_2n\bmod{r_2}\;\;\;\text{and}\;\;\;\;\ell_1c_2r_2\equiv -c_1c_2n\bmod{r_1}.
\end{align*}
Given $c_i$, $r_i$ and $n$, there are at most $(r_1,r_2)(c_1,r_2)(c_2,r_1)(n,r_1r_2)$ congruence classes $(a_1\bmod{r_1}, a_2\bmod{r_2})$ representing all possible solutions $(\ell_2,\ell_1)$. The number of $(\ell_1,\ell_2)$ falling in any congruence class is at most $O((1+L/p)^2)$. Hence the total contribution to $\Omega$ is bounded by
\begin{align*}
& \frac{\tilde{N}M^{3/2}}{m^2}\:\left(1+\frac{L^2}{p^2}\right)\:\mathop{\sum\sum}_{c_1,c_2\ll C} \:\mathop{\sum\sum}_{|r_1|,|r_2|\asymp p/m}\\
&\times  \sum_{0<|n|\ll p^2M^{1+\varepsilon}/\tilde{N}}\:(r_1,r_2)(c_1,r_2)(c_2,r_1)(n,r_1r_2).
\end{align*}
This is dominated by
\begin{align*}
& \frac{p^4M^{5/2}C^2}{m^4}\:\left(1+\frac{L^2}{p^2}\right).
\end{align*}
We conclude that
\begin{align*}
\Omega &\ll \frac{\tilde{N}M^{2}LpC}{m^2}+ \frac{\tilde{N}M^{3/2}L^2pC^2}{m^2}+\frac{p^4M^{5/2}C^2}{m^4}\:\left(1+\frac{L^2}{p^2}\right),
\end{align*}
from which it follows that
\begin{align*}
\mathcal{O}^\dagger\ll \frac{N^2\tilde{N}^{1/2}}{M^3p^3}\;\left[\frac{\tilde{N}^{1/2}M^{3/2}L^{1/2}p}{N^{1/2}}+ \frac{\tilde{N}^{1/2}M^{7/4}Lp^{3/2}}{N}+\frac{p^3M^{9/4}}{N}\:\left(1+\frac{L}{p}\right)\right].
\end{align*}
This reduces to
\begin{align*}
\mathcal{O}^\dagger\ll \frac{N^{1/2}M^{1/2}}{L^{1/2}}+ p^{1/2}M^{3/4}+\frac{N^{1/2}M^{1/4}}{L^{1/2}}\:\left(p+L\right).
\end{align*}
The optimal value $L=M^{1/4}$ is obtained by equating the first and the last term. The optimal value $p=M^{\varepsilon} NL/M=NM^{-3/4+\varepsilon}$ is obtained by choosing the minimal possible value for $p$ such that the condition required to make the initial off-diagonal negligibly small is satisfied. It follows that $\mathcal{O}^\dagger\ll N^{1/2}M^{3/8+\varepsilon}$. The same bound holds for the dual off-diagonal $\mathcal{O}^\star$. Consequently, as $\min\{p,L\}\gg NM^{-3/4-\varepsilon}$, we get
\begin{align*}
S(N)\ll N^{1/2}M^{3/8+\varepsilon} +N^{3/2}M^{-9/8+\varepsilon}+M^{3/4+\varepsilon}.
\end{align*}
This proves the first theorem, when we take $\theta=3/4$. Now observe that we have not used the cuspidality of $g$. If we replace $\lambda_g(n)$ by $d(n)$, the analysis goes through - only one has to replace the Rankin-Selberg $L$-function at the beginning of Section~\ref{fe} by $L(s,g)^2$. Thus we get the second theorem.

\bigskip




\begin{thebibliography}{99}


\bibitem{Agar} K. Aggarwal: $t$-aspect subconvexity for $GL(2)$ $L$-functions. arXiv:1707.07027


\bibitem{B} V. Blomer:
Subconvexity for twisted $L$-functions on $GL(3)$. Amer. J. Math. \textbf{134} (2012), 1385--1421.

\bibitem{BH} V. Blomer; G. Harcos: Hybrid bounds for twisted $L$-functions. J. reine angew. Math. \textbf{621} (2008), 53--79

\bibitem{BHM} V. Blomer; G. Harcos; P. Michel: A Burgess-like subconvex bound for twisted L-functions (with Appendix 2 by Z. Mao), Forum Math. \textbf{19} (2007), 61--105.

\bibitem{BMi} V. Blomer; G. Milicevic: p-adic analytic twists and strong subconvexity
Annales scientifiques de l'ENS \textbf{48} (2015), 561--605

\bibitem{BM} T.D. Browning; R. Munshi:
Rational points on singular intersections of quadrics.
Compositio Math., \textbf{149} (2013), 1457--1494.

\bibitem{Bu} D.A. Burgess:
On character sums and $L$-series. II. Proc. London Math. Soc. \textbf{13} (1963), 524--536.



\bibitem{CI} B. Conrey; H. Iwaniec: The cubic moment of central values
of automorphic L-functions. Annals of Mathematics, \textbf{151} (2000), 1175--1216.



\bibitem{DFI-1} W. Duke; J.B. Friedlander; H. Iwaniec: 
Bounds for automorphic L-functions.
Invent. Math. \textbf{112} (1993), 1--8. 



\bibitem{Go} A. Good: The square mean of Dirichlet series associated with cusp forms, Mathematika \textbf{29} (1982), 278--295.


\bibitem{IK} H. Iwaniec; E. Kowalski: Analytic Number Theory. Amer. Math. Soc. Coll. Publ. \textbf{53}, American Mathematical Society, Providence, RI, 2004.

\bibitem{KMV} E. Kowalski; P. Michel; J. Vanderkam: Rankin-Selberg L-functions in the level aspect,     Duke Math. J. \textbf{114} (2002), 123--191.


\bibitem{L} X. Li:
Bounds for $GL(3)\times GL(2)$ $L$-functions and $GL(3)$ $L$-functions. Annals of Math.  \textbf{173} (2011), 301--336.

\bibitem{Mil} G. Milicevic: Sub-Weyl subconvexity for Dirichlet L-functions to prime power moduli, Compositio Mathematica \textbf{152} (2016) 825--875.








\bibitem{Mu} R. Munshi: The circle method and bounds for $L$-functions - I. Math. Annalen, \textbf{358} (2014), 389--401.

\bibitem{Mu4} R. Munshi: The circle method and bounds for $L$-functions - II. Subconvexity for twists of $GL(3)$ $L$-functions. American J. Math., \textbf{137} (2015), 791--812.

\bibitem{Mu0} R. Munshi: The circle method and bounds for $L$-functions - III. $t$-aspect subconvexity for $GL(3)$ $L$-functions. J. Amer. Math. Soc., \textbf{28}  (2015), 913--938.
  
  \bibitem{Mu5} R. Munshi: The circle method and bounds for $L$-functions - IV. Subconvexity for twists of $GL(3)$ $L$-functions. Annals of Math., \textbf{182} (2015), 617--672.
  
  \bibitem{Mu6} R. Munshi: Pairs of quadrics in 11 variables.
Compositio Math., \textbf{151} (2015) 1189--1214.

\bibitem{Mu16} R. Munshi: Twists of $GL(3)$ $L$-functions. arXiv:1604.08000


\bibitem{MS} R. Munshi; S. Singh: Weyl bound for $p$-power twist of $GL(2)$ $L$-functions. arXiv:1706.03985

\bibitem{Sin} S. Singh: $t$ aspect subconvexity bound for $GL(2)$ $L$-functions. arXiv:1706.04977


\end{thebibliography}
\end{document}